\documentclass{amsart}
\usepackage[cp1251]{inputenc}
\usepackage[english,russian]{babel}
\usepackage{amsmath}
\usepackage{amssymb}
\usepackage{amsfonts}

\newtheorem{theorem}
{Theorem}

\newtheorem{proposition}{Proposition}

\begin {document}

\thispagestyle {empty}

\title {Corrections and additions  to my article \cite {G}}

\author {Gorbatsevich V.V.}

\maketitle

In my article \cite {G} I have found some inaccuracies and errors. Here I will describe in detail how to correct some incorrect  statements from \cite {G}, and which statements will have to be reformulated in a weaker form. At the same time, I decided not to limit myself to just pointing out corrections, but to give a more or less coherent text. Some definitions and constructions used will be given, but for details see \cite {G}. Here are some new results that expand the content of the article \cite {G}.

Recall that $\bf D$ denotes the algebra of dual numbers, its elements have the form $x+\epsilon y$, where $x,y \in \bf R$, and $\epsilon ^2=0$.

Let $M$ be some $\bf D$-manifold (i.e. it has a non-degenerate $\bf D$-structure $\mathcal E$) of real dimension $2n$. By  $\mathcal E$ we also denote the linear operator of a non-degenerate $ \bf D$-structure that acts on the tangent spaces of the manifold $M$ and, thus, on the space $T(M)$ of the tangent bundle over the manifold $M$. We denote the set of smooth vector fields on $M$ by $\Phi(M)$. They form an infinite-dimensional Lie algebra. By $\Phi_\mathcal E (M)$ we denote a Lie subalgebra  in $\Phi(M)$ consisting of $\bf D$-smooth fields. Vector fields from $\Phi_{\mathcal E} (M)$ have the form $\Sigma f_i(z) \partial/\partial z_i $, where $f_i(z)$ are $\bf D$-smooth functions. It is clear that $\Phi_{\mathcal E} (M)$ is a $\bf D$-Lie algebra (infinite-dimensional), i.e. it,  in addition to the standard identities for Lie algebras, also satisfies the identities $\mathcal E [Z, W] = [\mathcal EZ,W] = [Z, \mathcal E W]$ for arbitrary vector fields $Z, W \in \Phi_{\mathcal E} (M)$ (in other words,the commutation operation in $\Phi_{\mathcal E} (m)$ is $\bf D$-bilinear).

Consider the linear operator $\mathcal E_x$ of our $\bf D$-structure acting in the tangent space $T(M) _x$ of the manifold $M$ at the point $x\in M$. The operators $\mathcal E_x$ set an almost dual structure on $M$ (recall that $\mathcal E^2=0$), which was discussed in detail in \cite {G}. For this linear operator, consider the linear subspaces $Ker \mathcal E_x$ and $Im \mathcal E_x$ --- its kernel and image. For non-degenerate structures (we consider here only such structures ), these two subspaces coincide and their dimension is equal to $n$ --- the half the dimension of the tangent space. Given at each point of the manifold $M$, these subspaces form a smooth distribution on $M$, which we will call canonical. It is, as noted in \cite {G}, integrable (this fact is well known) and generates a smooth foliation on the manifold $M$, which we will also call canonical.

Since the operator $\mathcal E$ acts on every tangent space of the manifold $M$, it naturally acts on the Lie algebra $\Phi (M)$, as well as, as it is easy to check, on its Lie subalgebra $\Phi_{\mathcal E} (M)$. Important for us is the following property of this operator --- the intersection of its kernel $Ker \mathcal E$  (or, what is the same, its image) with $\Phi_{\mathcal E} (M)$ is an Abelian ideal in $\Phi_{\mathcal E} (m)$. Check it. We have $[\mathcal E Z, \mathcal EW] = \mathcal E^2[Z, W] =0$. This means that the image (aka the kernel) of the operator $\mathcal E$, when intersecting with $\Phi_{\mathcal E} (M)$, is an Abelian Lie  subalgebra. Further, the commutator $[\mathcal E Z, W] = \mathcal E [Z,W]$ lies in the image of this operator, and therefore this image is an ideal in $\Phi_{\mathcal E} (M)$.

Let us now proceed to the consideration of invariant $\bf D$-structures on manifolds. Let a certain Lie group $G$ acts on a $\bf D$-manifold $M$ (not necessarily having a $\bf D$-structure, but necessarily preserving such a structure on $M$). We will always assume that the action of the Lie group we are considering is locally effective. Therefore, it gives an embedding of the Lie algebra $g$ of the Lie group $G$ into the Lie algebra $\Phi_{\mathcal E} (M)$. We will also denote the image of this embedding by $g$. Since the action of a Lie group by our condition preserves the $\bf D$-structure on $M$ (thus the action is commutate with the operator $\mathcal E$), this image consists of left-invariant vector fields belonging to $\Phi_{\mathcal E}(M)$. The stationary subalgebra of this action at the point $m \in M$ consists of those vector fields of the Lie algebra $g$ that equals to 0 at the point $m$. If we consider a transitive action (we will consider such actions below), then all stationary subalgebras of this action are conjugate (with respect to the corresponding action on them of the Lie group $G$).  It is clear that $G$-invariant almost dual structure on $G/H$ is given by the linear operator $\mathcal E = \mathcal E(m_0)$, which acts on the tangent space to some fixed point $m_0$, as which we can take, for example, the point $eH$ if the homogeneous space $M$ is written as $G/H$ (where $H$ is a stationary subgroup of the point $m_0$ in $G$). The tangent space $T_{m_0}(M)$ to $M$ at the point $m_0$ we can identify with the quotient space $g/h$,  where $g,h$ are the Lie algebras of the Lie groups $G,H$ respectively. For an invariant almost dual structure, the canonical distribution is, of course, an invariant distribution. Also, the corresponding canonical foliation on a manifold is $G$-invariant one, i.e. the action of elements of the Lie group $G$ permutates the leaves of this foliation among themselves.  Note that leafs of canonical foliation  are connected submanifolds in $M$ of half dimension, but they can be nonclosed.

For a canonical foliation of a dual manifold, commutation of arbitrary vector fields on $M$ with fields tangent to the leafs of the canonical foliation will be tangent to this foliation.  In other words, the Lie algebra of vector fields tangent to the canonical foliation is an ideal in the Lie algebra of all vector fields on $M$. From the fact that the Lie algebra $\mathcal E(\Phi(M))$ is an Abelian ideal in $\Phi(M)$, it follows, in particular, that the commutator of two vector fields tangent to the canonical foliation leafs is zero.

Denote by $F$ a Lie subgroup in $G$ that keeps invariant the leaf of the canonical foliation  passing through the point $eH \in G/H$ . It is clear that $F$ is a Lie subgroup in $G$ (in general it is  not closed), and it contains a stationary Lie subgroup $H$. Due to the transitivity of the action of the Lie group $G$ on $M$ and the invariance of the canonical foliation, the natural action of the Lie group $F$ on the specified leaf of the foliation  is transitive, and therefore this leaf can be represented as $F/H$. From the definition of the Lie subgroup $F$, it is immediately clear that the corresponding Lie subalgebra $f \subset g$ contains the Lie subalgebra $g \cap Ker \mathcal E$, which is (as follows from the above) an Abelian ideal in the Lie algebra $g$ of the Lie group $G$.

Consider in more detail the stationary Lie subalgebra $h$ and the Lie algebra $f$ of the Lie group $F$. In \cite {G} (statement 2 of Proposition 4), it was stated that $f$ is an Abelian ideal in $g$. This statement is incorrect and it is the source of all incorrect statements in \cite {G}. Here is a simple example to prove this incorrectness.

Consider the $\bf D$-Lie group $G=SL_2({\bf D}) \cdot {\bf D} ^2$ --- the semi-direct product of $\bf D$-Lie group $SL_2(\bf D)$ and the Abelian ideal ${\bf D}^2$, which is given by a natural homomorphism of the Lie group $SL_2(\bf D)$ to the group of linear transformations of the dual plane ${\bf D}^2$. The ideal ${\bf D}^2$ can be written as a direct sum of two real planes ${\bf D}^2= {\bf R}^2 +\epsilon {\bf R}^2$. As a stationary subgroup $H$ we take $\bf D$-subgroup $SL_2({\bf D}) \subset G$. It is clear that the homogeneous space $G/H$ is exactly ${\bf D}^2$. As it is easy to understand, in this case the subalgebra $f$ corresponding to the Lie subgroup $F$ (which is transitive on the canonical foliation leafs) has the form $sl_2 ({\bf D}) + \epsilon {\bf R}^2$ . And since the subalgebra $sl_2 (\bf D)$ is non-Abelian (its exact description see in \cite {G}), then we see that the Lie subalgebra $f$ in this example is non-Abelian. The following statement corrects and supplements statement 2 of Proposition 4 of \cite {G}.

\begin{proposition}
 (an extended analog of Proposition 4 from \cite {G})

Let $M=G/H$ be a dual homogeneous space of the Lie group $G$. Then

1. A stationary Lie subalgebra $h \subset g$ is an even-dimensional regular subalgebra, i.e. the restriction of the operator $\mathcal E$ on it has the highest possible rank (which equals to half of its dimension).

2. The Lie algebra $f$ (corresponding to the Lie subgroup $F$) is an ideal in $g$ and can be represented as an exact sum $f=h+p$ of two Lie subalgebras, where $p$ is some Abelian Lie subalgebra in $f$ and $p \cap h =0$.

3. Semisimple parts of the Lie algebra $f$ and its Lie subalgebra $h$ are isomorphic (and they can be considered as coincident).

4. For the radicals $r_f$ and $r_h$ of the Lie algebra $f$ and its Lie subalgebra $h$, the exact decomposition $r_f=p+r_h$ takes place.

5. Let $k$ be a maximal semisimple compact Lie subalgebra in the Lie algebra $g$, considered as a real Lie algebra. Then the intersection of $k \cap h$ with a stationary Lie subalgebra $h$ will be a semisimple compact ideal in $g$.\end{proposition}

Proof.

The proof of statement 1 of this Sentence is contained in \cite {G} and does not need to be corrected. But statement 2 corrects the incorrect statement 2 from \cite {G} and we go to its proof.

Consider some subspace $p \subset f$, which is a complement to $h$. When projecting $g \to T_{eH}(M)$, the subalgebra $f$ passes into $Ker \mathcal E_{eH}$. So we can choose $p$ to lie in $f \cap Ker \mathcal E$. But for such Lie subalgebra $p$ it will be $[p,p]= 0$ (because the Lie subalgebra $Ker \mathcal E$ is Abelian), i.e. $p$ is an Abelian Lie subalgebra in $f$, and is complement to the Lie subalgebra $h$. So we prove statement 2 (which replaces the incorrect statement $[f,f] \subset h$ from \cite {G}).

Statements 3,4,5 are new (they were absent in \cite {G}). Let's move on to their proofs.

Let's start with statement 3. We use the fact, noted above, that the Lie subalgebra $f$ contains the Lie subalgebra $g \cap Ker\mathcal E$, which is an Abelian ideal in $g$, and therefore in $f$. Consider the natural epimorphism $\alpha : f \to f / g \cap Ker \mathcal E$. Since $f=p+h$ and $p \subset Ker \mathcal E$, the epimorphism $\alpha$ must be epimorphic under the restriction on $h$. Therefore, given the Abelian property of the ideal $Ker \mathcal E$, we get that the semisimple parts of the Lie algebras $f$ and $h$ must be isomorphic to each other (and the images of restrictions on them of the epimorphism $\alpha$ will be conjugate). But then the semisimple parts of the Lie subalgebras $f,h$ are also conjugate, and therefore they can be considered as coincident. This proves statement 3.

We proceed now to the statement 4. Consider the Levi decompositions of the Lie subalgebras $f$ and $h$. Since the semisimple parts of these Lie subalgebras coincide by virtue of the already proven statement 3, then these decompositions have the forms $f=s+r_f, h=s+r_h$, where $s$ --- their common semisimple part (the Levy factor), and $r_f, r_h$ --- their radicals.

Note that the Abelian ideal $g \cap Ker \mathcal E$ lies in $f$, and therefore the Lie subalgebra $p$ contained also in $r_f$. We show that the radical $r_h$ of the Lie subalgebra $h$ also lies in $r_f$. Assume that this is not the case and consider the natural epimorphism $\beta : f \to f/r_f =s$ of the Lie algebra $f$ onto its semisimple part. It follows from the above that the restriction of $ \beta$ on the Lie subalgebra $h$ will also be an epimorphism (due to the coincidence of semisimple parts of the Lie algebra $f$ and its Lie subalgebra $h$). If $r_h$ is not contained entirely in $r_f$, then $\beta (r_h)$ will be a nontrivial solvable ideal in a semisimple Lie algebra $s=f/r$ , which is impossible. The resulting contradiction shows that $r_h \subset r_f$.

Now using the facts that $f=p+h$, that semisimple parts of $f$ and $h$ are equal and that the subalgebras  $p$ and $r_h$ contained in $r_f$, it follows that $r_f=p+r_h$. It is clear that this is an exact decomposition. This completes the proof of statement 4. 

Statement 5 follows from the fact that $k \cap f$ is an ideal (automatically semisimple and compact) in $f$ and that the semisimple parts of the Lie subalgebras $f$ and $h$ can be considered as coincident by virtue of statement 4. This completes the proof of the entire Theorem 3. $\square$

For further information, it would be useful to know whether the decomposition $f=p+h$ from statement 2 of Proposition 1 is globalizable, i.e. whether there is a decomposition $F=P \cdot H$ of a Lie subgroup $F$ into the product of two its Lie subgroups (here $P$ is the connected Lie subgroup in $G$ corresponding to an Abelian Lie subalgebra $p$). From statements 3,4 of Proposition 1, it is clear that this will be true if and only if the exact decomposition $r_f=p+r_h$ of the radical of the Lie subalgebra $f$ into the exact sum of the Abelian lie subalgebra $p$ and the solvable Lie algebra $r_f$ is globalizable. It is known that in general, the decomposition of a solvable Lie algebra is not always globalizable. Apparently, this is also the case in our case. However, it is known that decompositions of nilpotent Lie algebras are always globalizable (see \cite {M}).

Let's consider the possibility of some reversal of statement 2 from Proposition 1. Assume that the Lie subalgebra $f$ is Abelian (which was incorrectly treated by \cite {G} as always fulfilled). Then all that is written on this subject on pages 46--47 in \cite {G} is true.

Further, Proposition 5 of \cite {G}, as noted there, is known in a more general situation. Therefore, the shorter proof given in that article (but, alas, based on the above incorrect statement) can be considered superfluous.

Our next correction relates to Theorem 3 of \cite {G}. It was formulated as follows:

\begin{theorem} (= Theorem 3 from \cite {G})

Let $M=G/H$ be a simply connected dual homogeneous space of a simply connected Lie group $G$. Then the manifold $M$ can be represented as the space of a vector $k$-dimensional bundle over a $k$-dimensional simply connected Lie group $G^\ast=G/F_0$. In particular, the manifold $M$ is homotopy equivalent to some compact Lie group.
\end{theorem}

The proof given in \cite {G} also uses the (incorrect in general) fact that the Lie subalgebra $f$ is Abelian. If we take this property of the Lie algebra $f$ as an additional assumption, then the statement of this Theorem is preserved. However, in a slightly modified form, it is true in a much more general situation. Namely, almost all arguments of the proof from \cite {G} pass if the decomposition $f=p+h$ is globalizable (see above for this condition). However, it is proved here that there is not a vector bundle over a compact Lie group, but just a bundle with a fiber ${\bf R}^n$. The statement about the homotopy type of a manifold $M$ does not change from this substitution.

Next, \cite {G} considers a number of examples that are not affected by an incorrect statement about the Abelian property of the Lie algebra $f$.

Section 4 of the article \cite {G} is devoted to compact dual homogeneous spaces. The statement of Proposition 7 there, concerning the case of a semisimple fundamental group, is in general (for a non-Abelian Lie subalgebra $f$) apparently incorrect, although it is preserved under the additional assumption of Abelian property of Lie algebra $f$. Here we consider a statement concerning compact homogeneous spaces with a finite fundamental group (a special case of the statement of Proposition 7 from \cite {G}, but not requiring the assumption of Abelian property of Lie algebra $f$).

\begin{proposition}

For an arbitrary compact homogeneous dual space $M=G/H$ of positive dimension its fundamental group $\pi_1(M)$ is infinite.
\end{proposition}

Proof.

Let $M=G/H$ be some compact homogeneous dual space. Assume that its fundamental group is finite. Then its universal covering manifold is also a compact homogeneous dual space, but already simply connected. So we can assume that $M$ is simply connected. By virtue of the well-known Montgomery theorem, the maximal compact semi-simple subgroup $K$ of the transitive Lie group $G$ is also transitive on $M$. Then, by virtue of statement 5 from Proposition 2 (see above) the intersection of the Lie algebra $k$ of the Lie group $K$ with the stationary Lie subalgebra $h \subset g$ is an ideal in $h$. But this means that when the Lie group $K$ acts on $M$, the stationary subgroup can be considered to be finite (considering only locally effective transitive actions). Thus we obtain that the simply connected manifold $M$ is a compact semisimple Lie group. But $g \cap Ker \mathcal E$ is the nontrivial Abelian ideal in a semisimple Lie algebra $g$, which is impossible. The resulting contradiction proves the statement of Proposition 3.$\square$

Next, \cite {G} considers a natural bundle for compact homogeneous spaces. The following statement is formulated and "proved" (using the same incorrect statement about Abelian property of algebra$f$):

\begin{theorem} (=Theorem 3 of \cite {G})

Let $M=G/H$ be a compact dual homogeneous space of a simply connected Lie group $G$. Then for some homogeneous space $M^\prime = G/H^\prime$, that covers $M$ finitely, the natural bundle is principal (with a compact semisimple structure group $K^\sharp$).
\end{theorem}

The proof given in \cite {G} can be slightly modified and it becomes correct. To do this, we need additionally to use statement 5 from Proposition 1, proven above.

Further, we note, that the argument given at the end of section 4 of the article \cite {G} about dual solvmanifolds, requires an additional assumption about the Abelian property of the Lie algebra $f$.

In the final section of the article \cite {G}, the properties of homogeneous $\bf D$-manifolds of dimensions 2,4 and 6 were considered. Also there were results on the structure of Lie groups transitive on such manifolds. Propositions 8,9 and 10 were based on a statement (incorrect) about the Abelian property of the Lie subalgebra $f$. Therefore, the descriptions given there of these transitive Lie groups are valid only under the additional assumption of Abelian property of Lie algebra $f$. In general, there are also Lie groups that are not included in the lists given in Propositions 8,9 and 10. The author intends to correct these lists in the future. For dimension 2 in the Proposition 8, except for claims (incomplete, as above) on the structure of transitive Lie groups, there is the following statement 2 about the topological description of homogeneous surfaces with $\bf D$-structure:

"2. The manifold  $M$ is diffeomorphic to the plane ${\bf R}^2$ or to the cylinder $s^2=S^1\times \bf R^1$ or to the torus $T^2$. On each of these three manifolds, structures of dual homogeneous spaces do exist."

This statement remains valid without requiring any adjustment. Also, comments in this section, that are not included in Propositions 8,9,10, do not need to be corrected.

\end{document}